\documentclass{article}

\usepackage[all]{xy}
\usepackage{amsmath,amssymb,amsthm}
\usepackage{stmaryrd}
\usepackage{authblk}

\usepackage{xcolor}

\newcommand{\jslattice}{\mathcal{S}^{\vee}_{0,1}}

\newtheorem{coro}{Corollary}
\newtheorem{defi}{Definition}

\newtheorem{rem}{Remark}
\newtheorem{prop}{Proposition}
\newtheorem{lem}{Lemma}

\newtheorem{theo}{Theorem}

\begin{document}

\title {Coextensive varieties via Central Elements}
\author[1,2]{W.J. Zuluaga Botero}
\affil[1]{\small Laboratoire J.A. Dieudonn\'e, CNRS, Universit\'e C\^ote d'Azur, 0618 Cedex 02, Nice--France}
\affil[2]{\small Departamento de Matem\'atica, Facultad de Ciencias Exactas, Universidad Nacional del Centro, Pinto 399, 7000 Tandil--Argentina}

%\thanks{This project has been funded by the European Research Council (ERC) under the European Union's Horizon 2020 research and innovation program (grant agreement No.670624). }

%\address{Laboratoire J. A. Dieudonn\'{e}, Universit\'{e} C\^ote d'Azur, Nice, France \\ Departamento de Matem\'aticas, Universidad Nacional del Centro de la Provincia de Buenos Aires, Tandil, Argentina}

\maketitle
\begin{abstract}
In this paper we use the theory of central elements in order to provide a characterization for coextensive varieties. In particular, if the variety is of finite type, congruence-permutable and its class of directly indecomposable members is universal, then coextensivity is equivalent to be a variety of shells.

\end{abstract}

%\tableofcontents

% NOTE: it is good practice to \label all headings (and proclamations) immediately

\section{Introduction}\label{Introduction}

In the category $\mathsf{Set}$ of sets and functions, coproducts are good enough to work with, in the sense that they have a good interaction with pullbacks. This means that coproducts are disjoint, the inclusions are monic and through the pullback of any function from a set $Z$ to $X+Y$ one may obtain a unique coproduct decomposition of $Z$. This idea motivated Schanuel \cite{S1990} and Lawvere \cite{L1990} to introduce extensive categories as categories $\mathcal{C}$ with finite coproducts and pullbacks in which the canonical functor $+:\mathcal{C}/X \times \mathcal{C}/X \rightarrow \mathcal{C}/(X+Y)$ is an equivalence. Therefore, one may think extensive categories (phrasing \cite{CW1993}) as ``\emph{categories in which coproducts exist and are well behaved}". Classical examples of extensive categories are the category $\mathsf{Set}$ of sets and functions, the category $\mathsf{Top}$ of topological spaces and continuous functions and the category $(\mathsf{CRing})^{op}$, where $\mathsf{CRing}$ is the category of commutative rings with unit and homomorphisms.  A category is called \emph{coextensive} if its opposite is extensive. If $\mathcal{V}$ is a variety (of universal algebras) then we say that $\mathcal{V}$ is coextensive if as an algebraic category is coextensive.

Varieties in which factor congruences of every algebra (i.e. kernel congruences of product projections) form a distributive sublattice of the lattice of all congruences are called \emph{varieties with Boolean factor congruences}. These varieties where studied by Chang, J\'onsson and Tarski in \cite{CJT1964}, in relation with a more classical property in universal algebra called \emph{Strict refinement property}. This property roughly says that any isomorphism between a product of irreducible structures is uniquely determined by a family of isomorphisms between each factor. However, along the literature the evidence of varieties in which direct product decompositions are determined by certain elements of the algebras of the variety is abundant. For instance, one may consider the case of idempotent elements in rings with unit, complemented elements in bounded distributive lattices, Boolean elements in residuated lattices \cite{KO2000}, MV-algebras \cite{C1958}, BL-algebras \cite{H1998}, Heyting algebras, etc. This fact implies that in all these cases such elements concentrate all the information concerning direct product decompositions. Within the context of varieties in which the universal congruence of each algebra is compact, in \cite{V1996_2} Vaggione introduces the concept of \emph{central element}. In \cite{SV2009} it is proved that under the aforementioned condition, varieties with Boolean factor congruences are equivalent to varieties in which factor congruences are definable by a first order formula.   
\\

This paper is motivated by the fact that some algebraic categories studied by the author (see \cite{CMZ2016} and \cite{Z2016}) were coextensive varieties and the proof of this fact always seemed to be related precisely with those elements which concentrate all the information about direct product decompositions. This observation generated the following question: ?`Is it possible to find a characterization of coextensive varieties in terms of the abovementioned elements? In the literature the are several caracterizations for coextensive categories which are given in purely categorical terms (see \cite{CW1993}, \cite{CPR2001} or \cite{H2020}) or even in syntactical terms (see \cite{B2020}). Nevertheless, as far as we know, none of these involves explicitly the notion of central element. This made the unifying approach provided by the theory of central elements a suitable way to study this problem.
\\

In this work we employ the theory of central elements in order to obtain a characterization of coextensive varieties (regardless the cardinality of its type) as Pierce varieties \cite{V1996} satisfying some concrete conditions (Theorem \ref{characterization coextensive}). In particular, if the variety is a variety of finite type such that every subalgebra of any subdirectly irreducible member is directly indecomposable, it turns out that coextensive varieties are precisely the Pierce varieties (Corollary \ref{extensive finite type}). Furthermore, if the variety is a congruence-permutable variety such that the class of its directly indecomposable members is a universal class, then coextensive varieties are exactly the varieties in which the universal congruence of each algebra is compact; or equivalently, varieties of shells (Corollary \ref{permutable coextensive}). Finally, we apply the obtained results to study the coextensivity of some particular families of varieties.

%The reader is assumed to be familiar with some standard topos theory as presented in \cite{MM2012} and \cite{J2002}. For standard notions in universal algebra the reader may consult \cite{MMT1987}.

\section{Preliminaries}\label{Preliminaries}

%\subsection{Notation and basic results}\label{Notation and basic results}

Let $A$ be a set and $k$ be a natural number. We write $\vec{a}$ for the elements of $A^{k}$. If $\vec{a} \in A^{k}$ and $\vec{b} \in B^{k}$, we write $[\vec{a}, \vec{b}]$ for the k-tuple $((a_{1} , b_{1} ), ..., (a_{k} , b_{k})) \in (A \times B)^{k}$. If $f:A\rightarrow B$ is a function and $\vec{a}\in A^{k}$, then we write $f(\vec{a})$ for the element $(f(a_{1}),...,f(a_{n}))\in B^{k}$. If $\mathbf{A}$ and $\mathbf{B}$ are algebras of the same type, we write $\mathbf{A}\leq \mathbf{B}$ to denote that $\mathbf{A}$ is subalgebra of $\mathbf{B}$. Let $\{\mathbf{A}_{j}\}_{j\in I}$ be a family of algebras of a given type. If there is no place to confusion, we write $\prod \mathbf{A}_{j}$ in place of $\prod_{j\in I} \mathbf{A}_{j}$. If $a,b\in \prod \mathbf{A}_{j}$, we write $E(a,b)$ for the set $\{j\in I\mid a(j)=b(j)\}$. If $\vec{a},\vec{b}\in (\prod \mathbf{A}_{j})^{k}$ then we write $E(\vec{a},\vec{b})$ for $\bigcap_{i=1}^{k}E(a_{j},b_{j})$. If $\mathbf{A}$ is an algebra of a given type and $X\subseteq A$ we write $\mathsf{Sg}^{\mathbf{A}}(X)$ for the subalgebra of $\mathbf{A}$ generated by $X$. We denote the congruence lattice of $\mathbf{A}$ by $\mathsf{Con}(\mathbf{A})$. If $\theta \in \mathsf{Con}(\mathbf{A})$, and $\vec{a}\in A^{k}$ we write $\vec{a}/\theta$ for the k-tuple $(a_{1}/\theta , ..., a_{k}/\theta)\in (A/\theta)^{k}$. The universal congruence on $\mathbf{A}$ is denoted by $\nabla^{\mathbf{A}}$ and $\Delta^{\mathbf{A}}$ denotes the identity congruence on $\mathbf{A}$ (or simply $\nabla$ and $\Delta$ when the context is clear). If $S\subseteq A\times A$, we write $\mathsf{Cg}^{\mathbf{A}}(S)$ for the least congruence generated by $S$. If $S=[\vec{a},\vec{b}]$,
we write $\mathsf{Cg}^{\mathbf{A}}(\vec{a},\vec{b})$ instead of $\mathsf{Cg}^{\mathbf{A}}([\vec{a},\vec{b}])$. We use $\mathsf{FC}(\mathbf{A})$ to denote the set of factor congruences of $\mathbf{A}$. 

\begin{lem}\label{embeddings are esencially inclusions}
Let $\mathbf{A}$ and $\mathbf{B}$ be algebras of the same type and let $e:\mathbf{A}\rightarrow \mathbf{B}$ be an embedding. Then, there exists an algebra $\mathbf{C}$ such that $\mathbf{A}\leq \mathbf{C}$ and $\mathbf{B}\cong \mathbf{C}$.
\end{lem}
\begin{proof}
Folklore.
\end{proof}

\begin{lem}\label{principal congruences homomorphisms}
Let $\mathbf{A}$ and $\mathbf{B}$ any algebras and let $f:\mathbf{A} \rightarrow \mathbf{B}$ be a homomorphism. Then $(a,b)\in \mathsf{Cg}^{\mathbf{A}}(\vec{a},\vec{b})$ implies $(f(a),f(b))\in \mathsf{Cg}^{\mathbf{B}}(f(\vec{a}),f(\vec{b}))$.
\end{lem}
\begin{proof}
Folklore.
\end{proof}

Let $\mathbf{A}$ and $\mathbf{B}$ be algebras of the same type and $\vec{a},\vec{b}\in A^{k}$. We say that a homomorphism $f:\mathbf{A} \rightarrow \mathbf{B}$ has the \emph{universal property of identify $\vec{a}$ with $\vec{b}$}, if for every homomorphism $g:\mathbf{A}\rightarrow \mathbf{C}$ such that $g(a_{i})=g(b_{i})$ for every $1\leq i\leq k$, there exists a unique homomorphism $h:\mathbf{B}\rightarrow \mathbf{C}$, making the diagram

\begin{displaymath}
\xymatrix{
\mathbf{A}  \ar[r]^-{f} \ar[dr]_-{g} & \mathbf{B}  \ar@{-->}[d]^-{h}
\\
 & \mathbf{C}
}
\end{displaymath}
\noindent
commute.

\begin{lem}\label{universal property principal congruences} Let $\mathbf{A}$ be an algebra and $S\subseteq A$. Then, the canonical homomorphism $\nu:\mathbf{A}\rightarrow \mathbf{A}/\mathsf{Cg}^{\mathbf{A}}(\vec{a},\vec{b})$ has the universal property of identify $\vec{a}$ with $\vec{b}$. 
\end{lem}
\begin{proof}
Apply Lemma \ref{principal congruences homomorphisms}.
\end{proof}

\begin{lem}\label{corollary universal property} 
Let $k$ be a positive number, $\mathbf{A}$ and $\mathbf{B}$ be algebras of the same type,  $f:\mathbf{A}\rightarrow \mathbf{B}$ be a homomorphism and $\nu$ and $\mu$ be the canonical homomorphisms. Then, for every $\vec{a},\vec{b} \in A^{k}$, the diagram
\begin{displaymath}
\xymatrix{
\mathbf{A} \ar[r]^-{\nu} \ar[d]_-{f} & \mathbf{A}/\mathsf{Cg}^{\mathbf{A}}(\vec{a},\vec{b}) \ar[d]^-{\psi}
\\
\mathbf{B} \ar[r]_-{\mu} & \mathbf{B}/\mathsf{Cg}^{\mathbf{A}}(f(\vec{a}),f(\vec{b}))
}
\end{displaymath}
\noindent
in where $\psi$ is the homomorphism raising from Lemma \ref{principal congruences homomorphisms}, is a pushout. 
\end{lem}
\begin{proof}
Apply Lemma \ref{universal property principal congruences}.
\end{proof}

Given a class $\mathcal{K}$ of algebras, we use $\mathbb{S}(\mathcal{K})$, $\mathbb{P}_{u}(\mathcal{K})$ and $\mathbb{V}(\mathcal{K})$ to denote the class of subalgebras and ultraproducts of elements of $\mathcal{K}$ and the variety generated by $\mathcal{K}$. For a
variety $\mathcal{V}$, we write $\mathcal{V}_{SI}$ and $\mathcal{V}_{DI}$ to
denote the classes of subdirectly irreducible and directly indecomposable
members of $\mathcal{V}$. Given a variety $\mathcal{V}$ and a set $X$ of variables we use $\mathbf{F}_{\mathcal{V}}(X)$ for the free algebra
of $\mathcal{V}$ freely generated by $X$. In particular, if $X=\{x_{1},...,x_{m}\}$, then we write $\mathbf{F}_{\mathcal{V}}(x_{1},...,x_{m})$ instead of $\mathbf{F}_{\mathcal{V}}(\{x_{1},...,x_{m}\})$. We stress that in this paper, by variety we mean \emph{variety with at least a constant symbol}. 
\\

A variety $\mathcal{V}$ has the \emph{Fraser-Horn property} (FHP, for
short) \cite{FH1970} if for every $\mathbf{A}_{1},\mathbf{A}_{2}\in \mathcal{%
V}$, it is the case that every congruence $\theta $ in $\mathbf{A}_{1}\times 
\mathbf{A}_{2}$ is the product congruence $\theta _{1}\times \theta _{2}$
for some congruences $\theta _{1}$ of $\mathbf{A}_{1}$ and $\theta _{2}$ of $%
\mathbf{A}_{2}$.

In \cite{C1971}, Comer defines the \emph{Pierce sheaf} of an algebra $\mathbf{A}
$ with Boolean factor congruences to be sheaf space of algebras $(F, \pi, X)$%
, where $X$ is the Stone space of maximal ideals of the Boolean algebra of
factor congruences of $\mathbf{A}$, $F$ is the disjoint union of the sets $%
\mathbf{A}/\bigcup \mathbf{m}$, with $\mathbf{m}\in X$, the map $\pi$ is
defined by $\pi(a/\bigcup \mathbf{m})=\mathbf{m}$, and $F$ is endowed with
the topology generated by the sets $\{a/\mathbf{m} \colon \mathbf{m}\in U\}$%
, with $a\in A$ and $U$ belonging to the topology of $X$. The natural map $%
\mathbf{A}\rightarrow \Pi \{A/\mathbf{m}\colon \mathbf{m}\in X \}$ is an
isomorphism between $\mathbf{A}$ and $\Gamma (X,F)$. In what follows, we
will refer to the $\mathbf{A}/\bigcup \mathbf{m}$'s as the \emph{Pierce
stalks of} $\mathbf{A}$. \newline

Let $\mathcal{L}$ be a first order language. If a $\mathcal{L}$-formula $\varphi (\vec{x})$ has the form 
\begin{equation*}
\forall \vec{y} \bigwedge_{j=1}^{n}p_{j}(\vec{x},\vec{y})=q_{j}(\vec{x},\vec{y}),
\end{equation*}%
for some positive number $n$ and terms $p_{j}(\vec{x},\vec{y})$ and $q_{j}(\vec{x},\vec{y})$ in $\mathcal{L}$, then we say that $\varphi (\vec{x})$ is a \emph{($\forall \bigwedge p=q$)-formula}. In a similar manner we define ($\bigwedge p=q$)-formulas. 

Let $\mathcal{L}$ be a first order language and $\mathcal{K}$ be a class of $%
\mathcal{L}$-structures. If $R\in \mathcal{L}$ is a $n$-ary relation
symbol, we say that a formula $\varphi (x_{1},...,x_{n})$ \emph{defines} $R$ 
\emph{in} $\mathcal{K}$ if%
\begin{equation*}
\mathcal{K}\vDash \varphi (\vec{x})\leftrightarrow R(\vec{x})\text{.}
\end{equation*}
In particular, if a ($\bigwedge p=q$)-formula defines $R$, we say that $R$ is \emph{equationally definable}.

\section{Central elements}

A \textit{variety with} $\vec{0}$ \textit{and} $\vec{1}$ is a variety $%
\mathcal{V}$ in which there are 0-ary terms $0_{1},\ldots ,0_{N},$ $%
1_{1},\ldots ,1_{N}$ such that $\mathcal{V}\vDash \vec{0}=\vec{1}\rightarrow
x=y$, where $\vec{0}=(0_{1},\ldots ,0_{N})$ and $\vec{1}=(1_{1},\ldots
,1_{N})$. The terms $\vec{0}$ and $\vec{1}$ are analogue, in a rather
general manner, to identity (top) and null (bottom) elements in rings
(lattices), and its existence in a variety, when the language has at least a
constant symbol, is equivalent to the fact that no non-trivial algebra in
the variety has a trivial subalgebra (see \cite{CV2012}). 

If $\mathbf{A}\in \mathcal{V}$, then we say that $\vec{e}\in A^{N}$ is a 
\textit{central element} of $\mathbf{A}$ if there exists an isomorphism $%
\mathbf{A}\rightarrow \mathbf{A}_{1}\times \mathbf{A}_{2}$ such that $\vec{e}%
\mapsto \lbrack \vec{0},\vec{1}]$. Also, we say that $\vec{e},\vec{f}\in A^{N}$ are a \emph{pair of complementary central elements} of $\mathbf{A}$ if there exists an isomorphism $\tau: \mathbf{A}\rightarrow \mathbf{A}_{1}\times \mathbf{A}_{2}$ such that $\tau(e)=[\vec{0}, \vec{1}]$ and $\tau(f)=[\vec{1}, \vec{0}]$. We use $Z(\mathbf{A})$ to denote the
set of central elements of $\mathbf{A}$ and $\vec{e}\diamond_{\mathbf{A}}\vec{f}$ to denote that $\vec{e}$ and $\vec{f}$ are complementary central elements of $\mathbf{A}$. It is clear from the above definition that for every $\mathbf{A}\in \mathcal{V}$ the set $\{(\vec{e},\vec{f})\in A^{2N}\colon \vec{e}\diamond_{\mathbf{A}}\vec{f}\}$ defines a $2N$-ary relation on $\mathbf{A}$.

Central elements are a generalization of both central idempotent elements in
rings with identity and neutral complemented elements in bounded lattices.
In these classical cases it is well known that central elements
concentrate the information concerning the direct product representations. 

Now, consider the following property:

\begin{enumerate}
\item[(L)] There is a first order formula $\lambda(x,y,\vec{z})$ such that, for every $\mathbf{A},\mathbf{B}\in \mathcal{V}$, $a,b\in A$ and $c,d\in B$ 
\[\mathbf{A}\times \mathbf{B}\models \lambda((a,c), (b,d), [\vec{0},\vec{1}])\;\textrm{iff}\; a=b. \] 
%\item[(R)] There is a first order formula $\rho(x,y,\vec{z})$ such that, for every $\mathbf{A},\mathbf{B}\in \mathcal{V}$, $a,b\in A$ and $c,d\in B$ 
%\[\mathbf{A}\times \mathbf{B}\models \rho((a,c), (b,d),[\vec{0},\vec{1}])\;\textrm{iff}\; c=d. \]
%\item[(W)] There is a a first order formula $\omega(x,y,\vec{z},\vec{w})$ such that, for every $\mathbf{A},\mathbf{B}\in \mathcal{V}$, $a,b\in A$ and $c,d\in B$ 
%\[\mathbf{A}\times \mathbf{B}\models \lambda((a,c), (b,d),[\vec{0},\vec{1}],[\vec{1},\vec{0}])\;\textrm{iff}\; a=b. \]

\end{enumerate} 
We say that a variety $\mathcal{V}$ has
\textit{definable factor congruences} if $(L)$ holds. In Theorem 1.1 of \cite{SV2009} it is proved that varieties $\mathcal{V}$
with $\vec{0}$ and $\vec{1}$ in which central elements determine the direct
product representations are precisely the varieties with definable factor congruences or equivalently, the varieties with Boolean factor congruences, i.e. those varieties whose set of factor congruences
of any algebra of $\mathcal{V}$ is a distributive sublattice of its congruence lattice. Moreover, in the same Theorem it is also proved that for every $\mathbf{A}\in \mathcal{V}$, the map $\mathsf{FC}(\mathbf{A}) \rightarrow   Z(\mathbf{A})$ that assigns to every factor congruence $\theta$ the unique $\vec{e}\in A^{N}$ satisfying $\vec{e}\equiv \vec{0}(\theta )$ and $\vec{e}\equiv \vec{1}(\theta^{\ast})$\footnote{We write $\vec{a}\equiv \vec{b}(\theta )$ to express that $(a_{i},b_{i})\in
\theta ,i=1,...,N$.} (where $\theta ^{\ast }$ is the complement of $\theta $ in $\mathsf{FC}(\mathbf{A})$)
is bijective. For every $\vec{e}\in Z(\mathbf{A})$ we write $\theta^{\mathbf{A}}_{\vec{0},\vec{e}}$ for the corresponding factor congruence obtained through the latter bijection. The latter, enables to endow $ Z(\mathbf{A})$ with a Boolean algebra structure. We write $\mathbf{Z}(\mathbf{A})$ for this Boolean algebra.
%See \cite{BV2020} for a non constructive short proof of this fact.
\\

For the details of the proofs of Lemmas \ref{In FHP surjections give
homomorphisms} and \ref{Products Center FHP}, the reader may consult the
proofs of the items $(a)$ and $(b)$ of Lemma 4 in \cite{V1999}.

\begin{lem}
\label{In FHP surjections give homomorphisms} Let $\mathcal{V}$ be a variety
with $\vec{0}$ and $\vec{1}$ with the FHP and let $\mathbf{A},\mathbf{B}\in 
\mathcal{V}$. If $f:\mathbf{A}\rightarrow \mathbf{B}$ is a surjective
homomorphism, then the map $Z(\mathbf{A})\rightarrow Z(\mathbf{B})$ defined
by $\vec{e}\mapsto (f(e_{1}),...,f(e_{N}))$, is a homomorphism from $\mathbf{%
Z}(\mathbf{A})$ to $\mathbf{Z}(\mathbf{B})$.
\end{lem}

\begin{lem}
\label{Products Center FHP} Let $\mathcal{V}$ be a variety with $\vec{0}$
and $\vec{1}$ with the FHP, $\{\mathbf{A}_{i}\}_{i\in I}$ be a family of
members of $\mathcal{V}$ and $\vec{e} \in (\prod \mathbf{A}_{i})^{N}$. Then, 
$\vec{e}\in Z(\prod \mathbf{A}_{i})$ if and only if $\vec{e}(i)\in Z(\mathbf{%
A}_{i})$ for every $i\in I$. Moreover, $\mathbf{Z}(\prod \mathbf{A}_{i})$ is
naturally isomorphic to $\prod \mathbf{Z}(\mathbf{A}_{i})$.
\end{lem}

Let $\mathcal{V}$ be a variety with Boolean factor congruences. If $\mathbf{A},\mathbf{B}\in \mathcal{V}$ and $f:\mathbf{A}\rightarrow \mathbf{B}$ is a homomorphism, we say that $f$ \emph{preserves central elements} if the map $f:Z(\mathbf{A})\rightarrow Z(\mathbf{B})$ is well defined; that is to say, for every $\vec{e}\in Z(\mathbf{A}),$ it follows that $f(\vec{e})\in Z(\mathbf{B})$. We say that $f$ \emph{preserves complementary central elements} if preserves central elements and for every $\vec{e}_{1}, \vec{e}_{2}\in Z(\mathbf{A})$, 
\[\vec{e}_{1}\diamond_{\mathbf{A}}\vec{e}_{2} \Rightarrow f(\vec{e}_{1})\diamond_{\mathbf{B}}f(\vec{e}_{2}). \]

\begin{defi}\label{Stability and C-Stability}
Let $\mathcal{V}$ be a variety with Boolean factor congruences. We say that $\mathcal{V}$ is stable by complements if every homomorphism preserves complementary central elements. 
\end{defi}

Observe that, since complements of central elements are unique, it follows that there is a bijection between $Z(\mathbf{A})$ and $\mathcal{K}_{\mathbf{A}}=\{(\vec{e},\vec{f})\in A^{2N}\mid \vec{e}\diamond_{\mathbf{A}}\vec{f}\}$. Then, if $f:\mathbf{A}\rightarrow \mathbf{B}$ is a homomorphism which preserves complementary central elements, it is clear that it also must preserve central elements. Classical examples of varieties with Boolean factor congruences in which every homomorphism preserves complementary central elements are the varieties ${\mathcal{R}}$ of commutative rings with unit and ${\mathcal{L}_{0,1}}$ of bounded distributive lattices. 
\\

%\begin{rem}\label{SC and CSC are not trivial}
We conclude this section by warning that Definition \ref{Stability and C-Stability} has to be taken carefully because there exist varieties with Boolean factor congruences such that their homomorphisms:
\begin{itemize}
\item[(1)] preserve central elements but does not preserve complementary central elements,
\item[(2)] does not preserve nor central elements nor complementary central elements. 
\end{itemize}

In order to illustrate the first situation, let $\jslattice$ be the variety of bounded join semilattices. Since $\jslattice$ is a variety with $0$ and $1$, and the formula $\varphi(x,y,z)= (x\vee z\approx y\vee z)$ satisfies condition (L) of the Introduction, it follows that $\jslattice$ is a variety with Boolean factor congruences. Let us consider the algebras $\mathbf{A}=\mathbf{2}\times \mathbf{2}$ and $\mathbf{B}=\mathbf{2}\times \mathbf{2}\times \mathbf{2}$ (with $\mathbf{2}$ the chain of two elements). Notice that $Z(\mathbf{A})=A$, $Z(\mathbf{B})=B$ and furthermore, $(1,0,0)\diamond_{\mathbf{B}}(0,1,1)$, $(0,1,0)\diamond_{\mathbf{B}}(1,0,1)$ and $(0,0,1)\diamond_{\mathbf{B}}(1,1,0)$. Let $\alpha:\mathbf{A}\rightarrow \mathbf{B}$ be the homomorphism defined by $\alpha(1,1)=(1,1,1)$, $\alpha(0,0)=(0,0,0)$, $\alpha(0,1)=(0,0,1)$ and $\alpha(1,0)=(0,1,1)$. It is clear that $\alpha$ preserves central elements but does not preserve complementary central elements. 

For the last situation, let $\mathcal{M}$ be the variety of bounded lattices. It is known (see \cite{V1999} and \cite{FH1970}) that $\mathcal{M}$ is a variety with Boolean factor congruences. If $\mathbf{C}=\mathbf{2}\times \mathbf{2}$ and $\mathbf{D}=\{0,1,a,b,c\}$, with $\{a,b,c\}$ non-comparable, it easily follows that $\textbf{C}$ is subalgebra of $\mathbf{D}$, but $\mathbf{C}$ is directly decomposable while $\mathbf{D}$ is not. So the inclusion does not preserve nor central elements nor complementary central elements. 
%\end{rem}

Hopefully, as we will see, in the particular case of the varieties we are dealing with, it will always be the case that homomorphisms preserve both complementary central elements and central elements (Lemma \ref{Equivalencia preservacion de morfismos}).

\section{Coextensive Varieties}\label{Coextensive Varieties}

We start by recalling that a category with finite products $\mathcal{C}$ is called coextensive if for each pair of objects $X, Y$ of $\mathcal{C}$ the canonical functor \[\times : X/\mathcal{C} \times Y/\mathcal{C} \rightarrow (X\times Y)/\mathcal{C}\] is an equivalence. The following result is the dual of Proposition 2.14 in \cite{CW1993}.

\begin{prop}\label{Coextensivity Proposition} 
A category $\mathcal{C}$ with finite products and pushouts along its projections is coextensive if and only if the following conditions hold:
\begin{enumerate}
\item (Products are codisjoint.) For every $X$ and $Y$, the projections $X\times Y \xrightarrow[]{p_{0}} X$ and $X\times Y\xrightarrow[]{p_{1}} Y$ are epic and the square below is a pushout
\begin{displaymath}
\xymatrix{
X\times Y \ar[r]^-{p_{1}} \ar[d]_-{p_{0}} & Y \ar[d]^-{!} 
\\
X \ar[r]_-{!}  & 1 
}
\end{displaymath}

\item (Products are pushout-stable.) For every $Y,X_{i}, Y_{i}$ with $i=0,1$ and $X_{0}\times X_{1} \xrightarrow[]{f} Y$, if the squares below are pushouts
\begin{displaymath}
\xymatrix{
X_{0} \ar[d]_-{h_{0}} & X_{0}\times X_{1} \ar[d]^-{f} \ar[l]_-{p_{0}} \ar[r]^-{p_{1}} & X_{1}  \ar[d]^-{h_{1}}
\\
Y_{0} & Y \ar[r]_-{y_{0}} \ar[l]^-{y_{1}} & Y_{1} 
}
\end{displaymath}
\noindent
then, the span $Y_{0}\leftarrow Y \rightarrow Y_{1}$ is a product.
\end{enumerate}
\end{prop}

\begin{rem}
Let $\mathcal{V}$ be a variety. In order to make the context clear, along this paper we will use $\mathcal{V}$ to denote both the variety and its associated algebraic category. 
\end{rem}

\begin{lem}\label{Codisjoint projections}
In every variety with $\vec{0}$ and $\vec{1}$, products are codisjoint.
\end{lem}
\begin{proof}
Let $\mathbf{A}\leftarrow \mathbf{A}\times \mathbf{B} \rightarrow \mathbf{B}$ be a product diagram in $\mathcal{V}$. Notice that projections send $[\vec{0},\vec{1}]\in \mathbf{A}\times \mathbf{B}$ to $\vec{0}$ in $\mathbf{A}$ and to $\vec{1}$ in $\mathbf{B}$, so $\vec{0}=\vec{1}$ in the pushout. Since $\mathcal{V}$ is a variety with $\vec{0}$ and $\vec{1}$, the result follows.
\end{proof}

We say that a variety $\mathcal{V}$ is a \textit{Pierce Variety} \cite{V1996} if there exist a positive natural number $N$, $0$-ary terms $0_{1}$,...,$0_{N}$, $1_{1}$,...,$1_{N}$ and a term $U(x,y,\vec{z},\vec{w})$, such that the following identities hold in $\mathcal{V}$:

\begin{displaymath}
\begin{array}{c}
U(x,y,\vec{0},\vec{1})=x \\
U(x,y,\vec{1},\vec{0})=y.
\end{array}
\end{displaymath}

In \cite{BV2016} it is shown that Pierce varieties play a crucial role  in the study of the definability of central elements and complementary central elements in some varieties with Boolean factor congruences in terms of some concrete set of formulas. A particular class of Pierce Varieties are \emph{varieties with a short decomposition term}; i.e. varieties with a term $u(x,y,\vec{z})$ and $0$-ary terms $0_{1}$,...,$0_{N}$, $1_{1}$,...,$1_{N}$ satisfying 
\begin{displaymath}
\begin{array}{c}
u(x,y,\vec{0})=x \\
u(x,y,\vec{1})=y.
\end{array}
\end{displaymath}
If $N=1$, varieties with a short decomposition term are called \emph{Church varieties} \cite{MS2008}. The study of Church varieties (c.f. \cite{MS2008},\cite{MS2010}, \cite{SLPK2013}) has been motivated by the fact that they comprise several mathematical structures arising from different fields of mathematics, which include, $\lambda$-abstraction algebras, combinatory algebras \cite{MS2010}, bounded integral commutative residuated lattices \cite{KO2000} and rings with unit, among others. In this context, the short decomposition term is called \emph{if-then-else term}.
\\

The following result provides some properties about Pierce varieties.

\begin{lem}\label{Properties of Pierce Varieties}
Let $\mathcal{V}$ be a Pierce variety.
\begin{enumerate}
\item $\mathcal{V}$ is a variety with $\vec{0}$ and $\vec{1}$.
\item $\mathcal{V}$ is a variety with Boolean factor congruences.
\item $\theta^{\mathbf{A}}_{\vec{0},\vec{e}}=\mathsf{Cg}^{\mathbf{A}}(\vec{0},\vec{e})$, for every $\mathbf{A}\in \mathcal{V}$ and every $\vec{e}\in Z(\mathbf{A})$.
\item $\mathcal{V}$ has the FHP.
\end{enumerate}
\end{lem}
\begin{proof}
For the proofs of (1), (2) and (3), see Propositions 3.1 and 3.2 in \cite{BV2016}. For (4), the reader may consult Theorem 5 in \cite{V1996}. 
\end{proof}

Let $\mathcal{V}$ be a Pierce variety and let $\mathcal{L}$ be the language of $\mathcal{V}$. Let $\mathbf{A}\in \mathcal{V}$ and $\vec{e},\vec{f}\in A^{N}$. Then from Theorem 5 of \cite{V1996} we can assure that $\vec{e}\diamond_{\mathbf{A}}\vec{f}$ if and only if, for every $a,b,c\in A$, $\vec{a},\vec{b}\in A^{m}$, $F\in \mathcal{L}$  and $1\leq i\leq N$ the following equations hold:
\begin{equation}
\begin{array}{l}
U(a,a,\vec{e},\vec{f}) =  a, 
\\
U(e_{i},1_{i},\vec{e},\vec{f}) = U(0_{i},e_{i},\vec{e},\vec{f}) = e_{i},
\\
U(1_{i},f_{i},\vec{e},\vec{f})  =  U(f_{i},0_{i},\vec{e},\vec{f})  =  f_{i},
\\
U(a,c,\vec{e},\vec{f})  =  U(a,U(b,c,\vec{e},\vec{f}),\vec{e},\vec{f})  =  U(U(a,b,\vec{e},\vec{f}),c ,\vec{e},\vec{f}),
\\
F(U(a_{1},b_{1},\vec{e},\vec{f}),...,U(a_{m},b_{m},\vec{e},\vec{f})) =  U(F(\vec{a}),F(\vec{b}),\vec{e},\vec{f}).
\end{array}
\label{Ecuaciones centrales}
\end{equation}

Hence it follows that in Pierce varieties, the relation $\vec{e}\diamond_{\mathbf{A}}\vec{f}$ is defined by a set of ($\forall \bigwedge p=q$)-formulas. An interesting outcome from the set of equations $(1)$ is given in Proposition 3.2 of \cite{BV2016}. Namely, the operations of $\mathbf{Z}(\mathbf{A})$ can be equationally described. In particular for every $\vec{e},\vec{f}\in Z(\mathbf{A})$, then $\vec{f}$ is the complement of $\vec{e}$ in $\mathbf{Z}(\mathbf{A})$ if and only if for every $1\leq i\leq N$: 
\begin{displaymath}
\begin{array}{l}
U(f_{i},1_{i},\vec{e},\vec{1})=U(f_{i},1_{i},\vec{1},\vec{e}),
\\
U(f_{i},0_{i},\vec{e},\vec{0})=U(f_{i},0_{i},\vec{0},\vec{e}).
\end{array}
\end{displaymath}

In the case of varieties with a short decomposition term, Proposition 3.3 of \cite{BV2016} brings an even simpler description for the complement of a central element. Namely, if $\vec{e},\vec{f}\in Z(\mathbf{A})$, then $\vec{f}$ is the complement of $\vec{e}$ in $\mathbf{Z}(\mathbf{A})$ if and only if for every $1\leq i\leq N$:   
\begin{displaymath}
\begin{array}{l}
f_{i}=u(1_{i},0_{i},\vec{e}).
\end{array}
\end{displaymath}

As a straightforward consequence from the latter, we obtain the following interesting property for Pierce varieties.

\begin{lem}\label{Equivalencia preservacion de morfismos}
Let $\mathcal{V}$ be a Pierce Variety and $\varphi:\mathbf{A}\rightarrow \mathbf{B}$ be a homomorphism. Then $\varphi$ preserves complementary central elements if and only if preserves central elements. 
\end{lem}

At this point one may be wandering if Pierce varieties are varieties with a short decomposition term. The answer to this question when regarding $\vec{0}$ and $\vec{1}$ fixed, is negative, as the following counterexample shows: Consider the variety $\mathcal{L}_{0,1}$ of bounded distributive lattices. It is clear that the term $U(x,y,z,w)=(x\vee z)\wedge (y\vee w)$, together with the constants $0$ and $1$, makes $\mathcal{L}_{0,1}$ a Pierce variety and suppose that $\mathcal{L}_{0,1}$ is a variety with a short decomposition term. Now take $\mathbf{A}$ as the distributive lattice of two elements, namely $0$ and $1$. Then, we have that for every $e,f\in Z(\mathbf{A})$, $f$ is the complement of $e$ in $\mathbf{Z}(\mathbf{A})$ if and only if $f=u^{\mathbf{A}}(1,0,e)$. So the equation $\neg x=u(1,0,x)$ holds in $\mathbf{A}$, where $\neg x$ denotes the complement of $x$. Then, since every distributive lattice $\mathbf{B}$ is a subdirect product with factors equal to $\mathbf{A}$, this would imply that $\mathbf{B}$ is a Boolean algebra, which is a contradiction.
\\

%\begin{lem}
%Let $\mathcal{V}$ be a variety with a short decomposition term. Then $\mathcal{V}$ is stable by complements if and only if the relation $\vec{e}\diamond_{\mathbf{A}}\vec{f}$ is equationally definable.
%\end{lem}
%\begin{proof}
%The ``if" part follows from item $(5)$ of Theorem 5.1 in \cite{BV2016}. The ``only if part" is immediate.
%\end{proof}

Now we are ready to show the main result of this paper.

\begin{theo}\label{characterization coextensive}
Let $\mathcal{V}$ a variety. Then, the following are equivalent:
\begin{enumerate}
\item $\mathcal{V}$ is coextensive.
\item $\mathcal{V}$ is a Pierce variety in which the relation $\vec{e}\diamond_{\mathbf{A}} \vec{f}$ is equationally definable.
\item $\mathcal{V}$ is a Pierce variety in which $\mathcal{V}_{DI}$ is a universal class.
\item $\mathcal{V}$ is a Pierce variety stable by complements.
\item $\mathcal{V}$ is variety with the FHP such that the stalks of the Pierce sheaf of every element of $\mathcal{V}$ are directly indecomposable.
\end{enumerate} 
\end{theo}
\begin{proof}
$(1)\Rightarrow (2)$. Let us assume $\mathcal{V}$ is coextensive and let $\mathcal{L}$ be the lenguage of $\mathcal{V}$. We start by proving that $\mathcal{V}$ is a Pierce Variety. To do so, notice that from the dual of Proposition 2.8 of \cite{CW1993}, the terminal object is strict. Since the terminal object is the trivial algebra, then in particular no nontrivial algebra of $\mathcal{V}$ has a trivial subalgebra. Therefore, from Proposition 2.3 of \cite{CV2012}, we get that there are unary terms $0_{1}(w)$,...,$0_{N}(w)$, $1_{1}(w)$,..., $1_{N}(w)$ such that 
\[\mathcal{V} \models \vec{0}(w)=\vec{1}(w) \rightarrow x=y, \]
where $x,y,w$ are distinct variables. Let $c\in \mathcal{L}$ be a constant symbol. Since $\mathcal{V} \models \vec{0}(c)=\vec{1}(c) \rightarrow x=y$, we can redefine $\vec{0}=\vec{0}(c)$ and $\vec{1}=\vec{1}(c)$.
\\

\noindent
Now, let $X=\{(x,x),(y,y),(0_{1},1_{1}),...,(0_{N},1_{N}),(1_{1},0_{1}),...,(1_{N},0_{N}) \}$, $Y=\{(0_{1},1_{1}),...,(0_{N},1_{N}),(1_{1},0_{1}),...,(1_{N},0_{N}) \}$, $\mathbf{A}=\mathbf{F}_{\mathcal{V}}(x,y)\times \mathbf{F}_{\mathcal{V}}(x,y)$, $\mathbf{B}=\mathsf{Sg}^{\mathbf{A}}(X)$ and $\mathbf{C}=\mathsf{Sg}^{\mathbf{A}}(Y)$. Observe that $\mathbf{C}\cong \mathbf{D}\times \mathbf{E}$, where \[\mathbf{D}=\mathsf{Sg}^{\mathbf{F}_{\mathcal{V}}(x,y)}(0_{1},...,0_{N},1_{1},...,1_{N})=\mathsf{Sg}^{\mathbf{F}_{\mathcal{V}}(x,y)}(1_{1},...,1_{N},0_{1},...,0_{N})=\mathbf{E}.\]
Consider the following diagram
\begin{displaymath}
\xymatrix{
\mathbf{D}\ar@{^{(}->}[d] \ar @{} [dr] |{(a)} & \ar[l]_-{\pi_{1}} \ar[r]^-{\pi_{2}} \ar@{^{(}->}[d] \mathbf{C} & \mathbf{E} \ar@{^{(}->}[d] \ar @{} [dl] |{(b)}
\\
\mathbf{F}_{\mathcal{V}}(x,y) \ar@{=}[d] & \mathbf{B} \ar[l]_-{p_{1}i} \ar[d]^-{i} \ar[r]^-{p_{2}i} & \mathbf{F}_{\mathcal{V}}(x,y) \ar@{=}[d] 
\\
\mathbf{F}_{\mathcal{V}}(x,y)  & \ar[l]^-{p_{1}} \ar[r]_-{p_{2}}  \mathbf{A} & \mathbf{F}_{\mathcal{V}}(x,y)
}
\end{displaymath} 
\noindent
where $\pi_{1},\pi_{2}, p_{1}, p_{2}$ are the respective projections, $i$ and each of the vertical arrows of the squares $(a)$ and $(b)$ are the respective inclusions, and the remaining verticals arrows are given by the identity. Because $\mathbf{B}$ contains $(x,x)$ and $(y,y)$ it is clear that $p_{1}i$ and $p_{2}i$ are epic. Now we prove that $(a)$ is a pushout. Since the commutativity of $(a)$ is clear, let $\mathbf{B}\xrightarrow{\alpha} \mathbf{H} \xleftarrow{\beta} \mathbf{D}$ be a cocone. Then, for every $\mathcal{L}$-term $t(x_{1},...,x_{2N})$, we obtain
\begin{displaymath}
\begin{array}{c}\label{ecuacion 1}
t^{\mathbf{D}}(\alpha(0_{1},1_{1}),...,\alpha(0_{N},1_{N}), \alpha(1_{1},0_{1}),...,\alpha(1_{N},0_{N})) = t^{\mathbf{D}}(0^{\mathbf{D}}_{1},...,0^{\mathbf{D}}_{N},1^{\mathbf{D}}_{1},...,1^{\mathbf{D}}_{N}),
\end{array}
\end{displaymath}
\noindent
so in particular we get that for every $1\leq j\leq N$,

\begin{eqnarray}\label{ecuacion 2}
\alpha(0_{j},1_{j})=0^{\mathbf{D}}_{j}, & & \alpha(1_{j},0_{j})=1^{\mathbf{D}}_{j}.
\end{eqnarray}
\noindent
Consider now the assignment $\{x,y\}\rightarrow D$ defined by $x \mapsto \alpha(x,x)$ and $y \mapsto \alpha(y,y)$ and let $\psi$ be the homomorphism which extends it. If we write $g$ for the inclusion of $\mathbf{D}$ in $\mathbf{F}_{\mathcal{V}}(x,y)$, it is easy to see that equation (\ref{ecuacion 2}) makes $\psi p_{1}i=\alpha$ and $\psi g=\beta$. The uniqueness of $\psi$ is granted from construction. Therefore, $(a)$ is a pushout. In a similar manner it can be proved that $(b)$ is a pushout. Hence, since $\mathcal{V}$ is coextensive, it is the case that the middle row of the diagram of above is a product diagram. Thus, since $(x,y)\in \mathbf{B}$, there exists an $\mathcal{L}$-term $U(x,y,\vec{z},\vec{w})$ such that 
\[(x,y)=U^{\mathbf{B}}((x,x),(y,y),(0_{1},1_{1}),...,(0_{N},1_{N}),...,(1_{1},0_{1}),(1_{N},0_{N})), \]
so we can conclude 
\begin{displaymath}
\begin{array}{ccc}
x=U(x,y,0_{1},...,0_{N},1_{1},...,1_{N}), &  & y=U(x,y,1_{1},...,1_{N},0_{1},...,0_{N}).
\end{array}
\end{displaymath}
\noindent
For the last part, first we prove that
\begin{itemize}
\item[(1)] There exists a set of $(\bigwedge p=q)$-formulas which defines $\vec{e}\diamond_{\mathbf{A}}\vec{f}$.
\end{itemize}

Let $\mathbf{A}\times \mathbf{B}\leq \mathbf{E}$, and let $i_{1}$ and $i_{2}$ be the pushouts of the inclusion $i:\mathbf{A}\times \mathbf{B}\rightarrow\mathbf{E}$ along $\pi_{\mathbf{A}}$ and $\pi_{\mathbf{B}}$, respectively.
\begin{displaymath}
\xymatrix{
\mathbf{A}\ar[d]_-{i_{1}} & \ar[l]_-{\pi_{\mathbf{A}}} \ar[r]^-{\pi_{\mathbf{B}}} \ar[d]^-{i} \mathbf{A}\times \mathbf{B} & \mathbf{B} \ar[d]^-{i_{2}}
\\
\mathbf{E}_{1}  & \ar[l]^-{\nu_{1}} \ar[r]_-{\nu_{2}}  \mathbf{E} & \mathbf{E}_{2} 
}
\end{displaymath} 

Since $\mathcal{V}$ is coextensive, then the inclusion coincides with $i_{1}\times i_{2}$ (so $i_{1}$ and $i_{2}$ are embeddings) and there exists an isomorphism $\tau:\mathbf{E}\rightarrow \mathbf{E}_{1} \times \mathbf{E}_{2}$. It is clear that $\tau i=i_{1}\times i_{2}$. Now, from Lemma \ref{embeddings are esencially inclusions}, there exist $\mathbf{A}_{1}$ and $\mathbf{B}_{1}$ such that $\mathbf{A}\leq \mathbf{A}_{1}$ and $\mathbf{E}_{1}\cong \mathbf{A}_{1}$; and such that $\mathbf{B}\leq \mathbf{B}_{1}$ and $\mathbf{B}_{1}\cong \mathbf{E}_{2}$. Therefore $\mathbf{E}_{1}\times \mathbf{E}_{2} \cong\mathbf{A}_{1}\times \mathbf{B}_{1}$ and if $\psi$ is such an isomorphism, it is clear that $\psi\tau i=\psi (i_{1}\times i_{2})$. Hence by (4) of Theorem 5.1 of \cite{BV2016}, (1) holds.
\\

Finally in order to prove that the set of (1) is finite, again by Theorem 5.1 of \cite{BV2016} it suffices to prove that $\mathbb{P}_{u}(\mathcal{V}_{SI})\subseteq \mathcal{V}_{DI}$. So let $\mathbf{B}=\prod \mathbf{A}_{j}/U\in \mathbb{P}_{u}(\mathcal{V}_{SI})$ and let $\vec{a}/U\in Z(\mathbf{B})$. Take $i_{0}\in I$ fixed and consider the diagram 

\begin{displaymath}
\xymatrix{
\mathbf{A}_{i_{0}} \ar[d]_-{f} & \ar[l]_-{\mu_{1}} \ar[r]^-{\mu_{2}} \ar[d]^-{\nu} \prod \mathbf{A}_{j} & \prod_{j\neq i_{0}} \mathbf{A}_{j} \ar[d]_-{g}
\\
\mathbf{B}/\mathsf{Cg}^{\mathbf{B}}(\vec{0}/U,\vec{a}/U)  & \ar[l]^-{p_{1}} \ar[r]_-{p_{2}}  \mathbf{B} & \mathbf{B}/\mathsf{Cg}^{\mathbf{B}}(\vec{1}/U,\vec{a}/U) 
}
\end{displaymath} 
Let $\vec{e}\in (\prod \mathbf{A}_{j})^{N}$ be the central element associated to the decomposition $\prod \mathbf{A}_{j}=A_{i_{0}}\times \prod_{j\neq i_{0}} \mathbf{A}_{j}$. From Lemma \ref{Properties of Pierce Varieties}, $\mathcal{V}$ has the FHP and since $\nu$ is surjective, then from Lemma \ref{In FHP surjections give homomorphisms},  $\nu(\vec{e})=\vec{e}/U$ belongs to $Z(\mathbf{B})$. Now, since $\mathcal{V}$ is coextensive, there exist morphisms $f$ and $g$ such that the squares of the diagram of above are pushouts. Then, from Lemma \ref{corollary universal property} it follows
\begin{displaymath}
\begin{array}{ccc}
\mathbf{B}/\mathsf{Cg}^{\mathbf{B}}(\vec{0}/U,\vec{a}/U)\cong \mathbf{B}/\mathsf{Cg}^{\mathbf{B}}(\vec{0}/U,\vec{e}/U) & \text{and} & \mathbf{B}/\mathsf{Cg}^{\mathbf{B}}(\vec{1}/U,\vec{a}/U)\cong \mathbf{B}/\mathsf{Cg}^{\mathbf{B}}(\vec{1}/U,\vec{e}/U),
\end{array}
\end{displaymath} 

which in turn implies that
\begin{displaymath}
\begin{array}{ccc}
\mathsf{Cg}^{\mathbf{B}}(\vec{0}/U,\vec{a}/U)=\mathsf{Cg}^{\mathbf{B}}(\vec{0}/U,\vec{e}/U) & \text{and} & \mathsf{Cg}^{\mathbf{B}}(\vec{1}/U,\vec{a}/U)=\mathsf{Cg}^{\mathbf{B}}(\vec{1}/U,\vec{e}/U).
\end{array}
\end{displaymath} 
Therefore by Corollary 4 of \cite{V1999}, $\vec{a}/U=\vec{e}/U$ and hence $E(\vec{a},\vec{e})\in U$. Moreover, from Lemma \ref{Products Center FHP} it follows $\vec{e}(j)\in \{\vec{0}(j),\vec{1}(j)\}$ for every $j\in I$. So $E(\vec{e},\vec{0})\cup E(\vec{e},\vec{1})=I\in U$. Then $E(\vec{e},\vec{0})\in U$ or $E(\vec{e},\vec{1})\in U$ because $U$ is an ultrafilter. For the first case, observe that $E(\vec{e},\vec{0})\cap E(\vec{a},\vec{e}) \subseteq E(\vec{a},\vec{0})$ so $E(\vec{a},\vec{0})\in U$, or equivalently $\vec{a}/U=\vec{0}/U$. In a similar fashion it can be proved that $\vec{a}/U=\vec{1}/U$. Hence, $\mathbf{B}\in \mathcal{V}_{DI}$ and from Theorem 5.1 of \cite{BV2016} the result follows.
\noindent

$(2)\Rightarrow (3)$. According to Theorem V.2.20 in \cite{BS1981}, for proving the statement, it is enough to check that $\mathcal{V}_{DI}$ is closed under the formation of isomorphisms, subalgebras and ultraproducts. Since the first two cases are easy, we only prove that $\mathbb{P}_{u}(\mathcal{V}_{DI})\subseteq \mathcal{V}_{DI}$. To do so, recall that we can take
\[\varphi(\vec{z},\vec{w})=\bigwedge_{i=1}^{n}p_{i}(\vec{z},\vec{w})=q_{i}(\vec{z},\vec{w}), \]
as a formula defining the relation $\vec{e}\diamond_{\mathbf{A}}\vec{f}$ in $\mathcal{V}$. So let $\{\mathbf{A}_{j}\}_{j\in I}$ be a family of directly indecomposable algebras of $\mathcal{V}$ and $U$ be an ultrafilter of $I$. Let $\mathbf{B}=\prod \mathbf{A}_{j} / U$ and suppose that $\mathbf{B}\notin \mathcal{V}_{DI}$, then there exist $\vec{e}/U,\vec{f}/U\in B^{N}$ such that $\vec{e}/U,\vec{f}/U\notin \{\vec{0}/U,\vec{1}/U\}$ and $\vec{e}/U\diamond_{\mathbf{B}}\vec{f}/U$. Hence, by \L o\'s 's  theorem:
\begin{displaymath}
\begin{array}{cccc}
\llbracket \varphi(\vec{e},\vec{f}) \rrbracket = \{j\in I\mid \mathbf{A}_{j}\models \varphi(\vec{e}(j),\vec{f}(j)) \} \in U.
\end{array}
\end{displaymath}  
Since $\mathbf{A}_{j}\in \mathcal{V}_{DI}$ for every $j\in I$, then it follows  
\[\llbracket \varphi(\vec{e},\vec{f}) \rrbracket \subseteq \{j\in I\mid \vec{e}(j),\vec{f}(j)\in \{\vec{0}(j),\vec{1}(j) \}\}=C\]
so due to $U$ is increasing, we have $C\in U$. Observe that $C=E\cup F$, where
\begin{displaymath}
\begin{array}{cccc}
E & = &\{j\in I\mid \vec{e}(j)=\vec{0}(j)\}\cap \{j\in I\mid \vec{f}(j)=\vec{1}(j)\},
\\
F & = &\{j\in I\mid \vec{e}(j)=\vec{1}(j)\}\cap \{j\in I\mid \vec{f}(j)=\vec{0}(j)\}.
\end{array}
\end{displaymath} 
Thus again, because $U$ is an ultrafilter we have $E\in U$ or $F\in U$. It is clear that this implies $\vec{e}/U,\vec{f}/U \in \{\vec{0}/U,\vec{1}/U\}$, which is absurd. 

$(3)\Rightarrow (4)$. Since $\mathcal{V}$ is a universal class we have $\mathbb{S}(\mathcal{V}_{SI})\subseteq \mathcal{V}_{DI}$ and $\mathbb{P}_{u}(\mathcal{V}_{SI})\subseteq \mathcal{V}_{DI}$, so from Theorem 5.1 of \cite{BV2016}, the relation $\vec{e}\diamond_{\mathbf{A}}\vec{f}$ is equationally definable. From latter it easy to see that $\mathcal{V}$ is stable by complements.

$(4)\Rightarrow (1)$. Notice that from Lemma \ref{Properties of Pierce Varieties}, $\mathcal{V}$ is a variety with $\vec{0}$ and $\vec{1}$. Then from Lemma \ref{Codisjoint projections} in $\mathcal{V}$, products are codisjoint.  We only remain to prove that in $\mathcal{V}$ products are pushout-stable. To do so, let $\mathbf{A},\mathbf{B}\in \mathcal{V}$ and $f:\mathbf{A}\rightarrow \mathbf{B}$ be a homomorphism. From Lemma \ref{Properties of Pierce Varieties}, $\mathcal{V}$ is a variety with Boolean factor congruences and for every $\vec{e}\in Z(\mathbf{A})$, $\theta^{\mathbf{A}}_{\vec{0},\vec{e}}=\mathsf{Cg}^{\mathbf{A}}(\vec{0},\vec{e})$. Therefore without any loss of generality we can assume $\vec{e}\diamond_{\mathbf{A}}\vec{f}$ and consider the following diagram in where $\mathbf{P}_{1}$ and $\mathbf{P}_{2}$ are the pushouts from the left and the right squares, respectively, and the $\mu_{i}$ are the canonical homomorphisms.

\begin{displaymath}
\xymatrix{
\mathbf{A}/\mathsf{Cg}^{\mathbf{A}}(\vec{0},\vec{e}) \ar[d] & \ar[l]_-{\mu_{1}} \ar[r]^-{\mu_{2}} \ar[d]^-{f} \mathbf{A} & \mathbf{A}/\mathsf{Cg}^{\mathbf{A}}(\vec{0},\vec{g}) \ar[d]
\\
\mathbf{P}_{1}  & \ar[l]^-{p_{1}} \ar[r]_-{p_{2}}  \mathbf{B} & \mathbf{P}_{2} 
}
\end{displaymath} 
\noindent
From assumption $f(\vec{e})\diamond_{\mathbf{B}} f(\vec{g})$ so $\mathsf{Cg}^{\mathbf{B}}(\vec{0},f(\vec{e}))$ and $\mathsf{Cg}^{\mathbf{B}}(\vec{0},f(\vec{g}))$ are complementary factor congruences of $\mathbf{B}$. By Lemma \ref{corollary universal property}, there exist isomorphisms $i_{1}:  \mathbf{B}/\mathsf{Cg}^{\mathbf{B}}(\vec{0},f(\vec{e})) \rightarrow \mathbf{P}_{1}$ and $i_{2}:\mathbf{B}/\mathsf{Cg}^{\mathbf{B}}(\vec{0},f(\vec{g})) \rightarrow \mathbf{P}_{2}$. This leads us to consider the following diagram:

\begin{displaymath}
\xymatrix{
\mathbf{A}/\mathsf{Cg}^{\mathbf{A}}(\vec{0},\vec{e}) \ar[d] & \ar[l]_-{\mu_{1}} \ar[r]^-{\mu_{2}} \ar[d]^-{f} \mathbf{A} & \mathbf{A}/\mathsf{Cg}^{\mathbf{A}}(\vec{0},\vec{g}) \ar[d]
\\
\mathbf{B}/\mathsf{Cg}^{\mathbf{B}}(\vec{0},f(\vec{e})) \ar[d]_-{i_{1}} & \mathbf{B} \ar[l]_-{\nu_{1}} \ar[d]^-{id_{\mathbf{B}}} \ar[r]^-{\nu_{2}} & \mathbf{B}/\mathsf{Cg}^{\mathbf{B}}(\vec{0},f(\vec{g})) \ar[d]^-{i_{2}} 
\\
\mathbf{P}_{1}  & \ar[l]^-{p_{1}} \ar[r]_-{p_{2}}  \mathbf{B} & \mathbf{P}_{2} 
}
\end{displaymath} 
\noindent

In where the $\nu_{i}$ are also the canonical homomorphisms. Observe that since the upper left and the outer left squares are pushouts then, the lower left square is a pushout, hence in particular, $p_{1}=i_{1}\nu_{1}$ and $p_{2}=i_{2}\nu_{2}$. Let $\mathbf{P}_{1}\xleftarrow{\alpha} \mathbf{C} \xrightarrow{\beta} \mathbf{P}_{2}$ be a span. Since the span $\mathbf{B}/\mathsf{Cg}^{\mathbf{B}}(\vec{0},f(\vec{e})) \xleftarrow{\nu_{1}}\mathbf{B} \xrightarrow{\nu_{2}} \mathbf{B}/\mathsf{Cg}^{\mathbf{B}}(\vec{0},f(\vec{g}))$ is a product, let $k=\langle i_{1}^{-1}\alpha, i_{2}^{-1}\beta\rangle$ be the homomorphism raising from the universal property of the product. It can be easily proved that $k$ is the only arrow such that $p_{1}k=\alpha$ and $p_{2}k=\beta$. I.e. the span $\mathbf{P}_{1}\xleftarrow{p_{1}} \mathbf{B} \xrightarrow{p_{2}} \mathbf{P}_{2}$ is a product. Hence, by Proposition \ref{Coextensivity Proposition}, the result holds.

$(3)\Leftrightarrow (5)$. This follows from Theorem 8 of \cite{V1996}.

\end{proof}

\begin{coro}\label{Locally finite}
Let $\mathcal{V}$ be a Pierce variety such that $\mathbb{S}(\mathcal{V}_{SI})\subseteq \mathcal{V}_{DI}$. If $\mathcal{V}$ is locally finite, then $\mathcal{V}$ is coextensive.
\end{coro}
\begin{proof}
From Corollary 5.2 of \cite{BV2016} we get that $\mathcal{V}_{DI}$ is a universal class. So by Theorem \ref{characterization coextensive} the result holds.
\end{proof}

It is no hard to see by hand that $\mathcal{L}_{0,1}$ is coextensive. Nevertheless, notice that as an application of Corollary \ref{Locally finite} we can obtain an alternative proof of this fact. It is well known that $\mathcal{L}_{0,1}$ is locally finite so since the only subdirectly irreducible member of this variety is the distributive lattice of two elements, from Corollary \ref{Locally finite} the result follows.  
\\

As we have saw at the beginning of this section, the varietes $\mathcal{S}_{0,1}^{\vee}$ of bounded join semilattices and $\mathcal{M}$ of bounded lattices are not stable by complements so from Theorem \ref{characterization coextensive} it follows that their are not coextensive varieties. 
\\

We stress that there are Pierce varieties in which the relation $\vec{e}\diamond_{\mathbf{A}}\vec{f}$ is not equationally definable, and therefore by Theorem \ref{characterization coextensive}, they are not coextensive. For each $k\geq 2$, let $C_{k}$ be the bounded chain of $k$ elements. Define $\mathbf{C}_{k}=(C_{k},\vee, \wedge, \Rightarrow, 0,1)$, where
\begin{displaymath}
x\Rightarrow y= \left\{ \begin{array}{ll}
             1 &  \text{if}\; x \leq y \\
             \\ 0 &  \text{otherwise.}
             \end{array}
   \right.
\end{displaymath}

Let $\mathcal{L}$ be the language of bounded distributive lattices expanded by adding a new binary function symbol $\Rightarrow$. Let $\mathcal{L}_{E}=\mathcal{L}\cup \{F_{2},F_{3}...\}$ where for each $n\geq 2$, $F_{n}$ is a $n$-ary relation symbol. For $k\geq 2$, let $\mathbf{G}_{k}=(\mathbf{C}_{k}\times \mathbf{C}_{k},f_{2},f_{3},...)$, where $f_{n}=F^{\mathbf{G}_{k}}_{n}$ and for each $n\neq k$
\begin{displaymath}
\begin{array}{cc}
f_{n}(x_{1},...,x_{n})=1, & \text{for every}\; x_{1},...,x_{n}\in C_{k}\times C_{k}
\end{array}
\end{displaymath}
and 
\begin{displaymath}
f_{k}(x_{1},...,x_{k})= \left\{ \begin{array}{ll}
             0 &  \text{if}\; x_{i}\neq x_{j}\; \text{for every}\; i\neq j, \\
             \\ 1 &  \text{otherwise.}
             \end{array}
   \right.
\end{displaymath}
  
Let $\mathcal{G}=\{\mathbf{G}_{k}\colon k\geq 2\}$. In \cite{BV2016} it is proved that $\mathbb{V}(\mathcal{G})$ is a variety with a short decomposition term given by $u(x,y,z)=(x\wedge \neg z)\vee (x\wedge z)$ in which there is no finite set of formulas defining the relation $\vec{e}\diamond_{\mathbf{A}}\vec{f}$. Hence, by Theorem \ref{characterization coextensive}, $\mathbb{V}(\mathcal{G})$ is not coextensive.

\begin{coro}\label{extensive finite type}
Let $\mathcal{V}$ be a variety of finite type such that $\mathbb{S}(\mathcal{V}_{SI})\subseteq \mathcal{V}_{DI}$. Then, $\mathcal{V}$ is a Pierce variety if and only if $\mathcal{V}$ is coextensive.
\end{coro}
\begin{proof}
If $\mathcal{V}$ is a coextensive variety the result follows from Theorem \ref{characterization coextensive}. On the other hand, let us assume that $\mathcal{V}$ is a Pierce variety. Since $\mathbb{S}(\mathcal{V}_{SI})\subseteq \mathcal{V}_{DI}$ from assumption, then from Theorem 5.1 of \cite{BV2016}, there exists a set $\Sigma$ of $(\bigwedge p=q)$-formulas which defines the relation $\vec{e}\diamond_{\mathbf{A}}\vec{f}$. But the type of $\mathcal{V}$ is finite, so $\Sigma$ is finite. Therefore, again from Theorem \ref{characterization coextensive}, $\mathcal{V}$ is coextensive.

\end{proof}

According to \cite{V1996_1} a \emph{variety of shells} is a variety $\mathcal{V}$ in which there is a positive number $N$, $0$-ary terms $0_{1}$,...,$0_{N}$, $1_{1}$,...,$1_{N}$, terms $f_{i}(x_{1},...,x_{n},y_{1},...,y_{n})$ and $g_{i}(x_{1},...,x_{n},y_{1},...,y_{n})$, with $1\leq i\leq N$ such that for every $i$:
\begin{displaymath}
\begin{array}{ll}
f_{i}(\vec{x},\vec{0})=f_{i}(\vec{0},\vec{x})=0_{i},
\\
f_{i}(\vec{x},\vec{1})=f_{i}(\vec{1},\vec{x})=x_{i},
\\
g_{i}(\vec{x},\vec{0})=g_{i}(\vec{0},\vec{x})=x_{i}.
\end{array}
\end{displaymath} 

It is usual to write the terms of above as $f_{i}(\vec{x},\vec{y})=\vec{x}\times_{i}\vec{y}$ and $g_{i}(\vec{x},\vec{y})=\vec{x}+_{i}\vec{y}$. With this notation, the latter equations adopt a more familiar form 
\begin{displaymath}
\begin{array}{ll}
\vec{x}\times_{i}\vec{0}=\vec{0}\times_{i}\vec{x}=0_{i},
\\
\vec{x}\times_{i}\vec{1}=\vec{1}\times_{i}\vec{x}=x_{i},
\\
\vec{x}+_{i}\vec{0}=\vec{0}+_{i}\vec{x}=x_{i}.
\end{array}
\end{displaymath}

We conclude this section with a characterization of a particular class of congruence-permutable varieties of shells.

\begin{coro}\label{permutable coextensive}
Let $\mathcal{V}$ be a congruence-permutable variety of finite type such that $\mathcal{V}_{DI}$ is a universal class. Then, the following are equivalent:
\begin{enumerate}
\item $\mathcal{V}$ is a variety of shells.
\item $\mathcal{V}$ is a variety with $\vec{0}$ and $\vec{1}$.
\item $\mathcal{V}$ is a Pierce variety.
\item $\mathcal{V}$ is a coextensive variety.
\end{enumerate}
\end{coro}
\begin{proof}
For (1)$\Leftrightarrow$(2)$\Leftrightarrow$(3) the reader may consult \cite{V1996_1}. Finally,   observe that since $\mathcal{V}_{DI}$ is a universal class, then $\mathbb{S}(\mathcal{V}_{SI})\subseteq \mathcal{V}_{DI}$, so (3)$\Leftrightarrow$(4) follows from Corollary \ref{extensive finite type}.
\end{proof}

\section{Applications}

In this section we use the techniques we have developed so far to provide some examples of coextensive varieties. In particular, we study the coextensivity of some classes of preprimal varieties. 
\\

A \emph{discriminator algebra} is a nontrivial algebra $\mathbf{A}$ for which there exists a term $t(x,y,z)$ in the language of $\mathbf{A}$ such that for every $a,b,c\in A$
\begin{displaymath}
\begin{array}{cccc}
t^{\mathbf{A}}(a,a,c)=c & \text{and} & t^{\mathbf{A}}(a,b,c)=a & \text{if}\; a\neq b.
\end{array}
\end{displaymath}
A \emph{discriminator variety} is a variety generated by a class of (similar) algebras
which are discriminator algebras with respect to the same term $t$.

\begin{prop}
Every discriminator variety is coextensive.
\end{prop}
\begin{proof}
In \cite{V2001} it was proved that the Pierce stalks of every member in a discriminator variety are directly indemposable. Then, the result follows from Theorem \ref{characterization coextensive}. 
\end{proof}

An algebra $\mathbf{P}$ is called \emph{primal} if $\mathbf{P}$ is finite and every finitary operation on $P$ is a term operation. A variety $\mathcal{V}$ is primal if it is generated by a primal algebra. It is well known that every primal variety is a discriminator variety, therefore we obtain the following result.  

\begin{coro}
Every primal variety is coextensive.
\end{coro}

An algebra $\mathbf{P}$ is called \emph{preprimal} if $\mathbf{P}
$ is finite and its clone $\mathrm{Clo}(\mathbf{P})$ is a
maximal clone. A \textit{preprimal variety} is a variety generated by a
preprimal algebra.

In \cite{R1970} Rosenberg characterizes the
maximal clones over a a finite base set in terms of relations. He described for each
preprimal algebra $\mathbf{P}$ an $m$-ary relation $\sigma$ on $%
P$ in such a way that the $n$-ary term-operations of $\mathbf{P}$ are precisely the $n
$-ary functions $f$ on $P$ which preserve $\sigma$. These relations belong to
one of the following seven types:

\begin{enumerate}
\item Permutations with cycles all the same prime length,

\item Proper subsets,

\item Prime-affine relations,

\item Bounded partial orders,

\item $h$-adic relations,

\item $h$-ary central relations, with $h\geq 2$,

\item Proper, non-trivial equivalence relations.
\end{enumerate}

If $\sigma$ is a relation of any of the aforementioned types we denote by $\mathbf{P}_{\sigma} $ the preprimal algebra whose universe is $P$ and whose fundamental
operations are all the functions preserving $\sigma$. We write $\mathbb{V}(\mathbf{P}_{\sigma })$ to denote the variety generated by $\mathbf{P}_{\sigma }$.

\begin{prop}
Any of the preprimal varieties of types 1., 2.,  and 7. are coextensive. Moreover, the preprimal varieties of the type 6. only are coextensive if $h=2$.
\end{prop}
\begin{proof}
In \cite{K2012} it was proved that $\mathbb{V}(\mathbf{P}_{\sigma })$ has the FHP for any $\sigma$ of the types mentioned in the statement. Furthermore, in the same paper it was proved that the Pierce stalks of every member of $\mathbb{V}(\mathbf{P}_{\sigma })$ with $\sigma$ of type 1. and 2. are directly indecomposable. In \cite{VZ2020} the same result was proved for $\mathbb{V}(\mathbf{P}_{\sigma })$ with $\sigma$ of type 7. and 6. with $h=2$. In the same paper it was also proved that for the case of $\mathbb{V}(\mathbf{P}_{\sigma })$ with $\sigma$ of type 6. with $h\geq 3$, there are Pierce stalks of $\mathbb{V}(\mathbf{P}_{\sigma })$ which are not directly indecomposable. Then from Theorem \ref{characterization coextensive} the result holds. 
\end{proof}

\section*{Acknowledgements}
This project has been funded by the European Research Council (ERC) under the European Union's Horizon 2020 research and innovation program (grant agreement No.670624) and by the CONICET [PIP 112-201501-00412].

\end{document}